\documentclass[11pt]{article}
\usepackage{graphicx,epsfig,amssymb,cite}

\date{}

\begin{document}
\title{Finite-time stability and stabilization  of  linear discrete
time-varying stochastic
systems  }
\author{ Tianliang~Zhang ${}^1$,
 Feiqi~Deng ${}^1$ \thanks{
 Corresponding author. Email: aufqdeng@scut.edu.cn},
 and Weihai~Zhang ${}^2$ \\
${}^1$ {\small School of Automation Science and Engineering,} \\
 {\small South China University of Technology,
Guangzhou 510640,   P. R. China}\\
${}^2$ {\small College of Information and Electrical Engineering,} \\
{\small Shandong University of Science and Technology,}  {\small
Qingdao  266510, P. R. China}\\
 }\maketitle

{\bf Abstract-} {This paper studies the finite-time stability and
stabilization of linear discrete time-varying stochastic systems with
multiplicative noise. Firstly, necessary  and  sufficient conditions
for  finite-time stability  are presented via state transition
matrix approach. Secondly,  this paper also develops the Lyapunov
function method   to study  finite-time stability and stabilization
of discrete time-varying stochastic systems based on matrix
inequalities and linear matrix inequalities (LMIs),  so as to Matlab
LMI Toolbox can be used. Two numerical examples are given to
illustrate the effectiveness of the proposed results. }


\section {Introduction}
As it is well-known that stability is the first
consideration in system analysis and synthesis. Since A. M. Lyapunov
published his  classical work \cite{Lyapunov} on  stability of
ordinary differential equations (ODEs) in 1892, Lyapunov's stability
theory has been one of the most important issues in mathematics and
modern control theory. In particular, Lyapunov's second method has
been extended to continuous-time stochastic It\^o-type differential
equations, and we refer the reader to \cite{hudeng,hasmin_80,mao}.
Stability in Lyapunov sense describes the asymptotic behaviour of
the  state trajectory as time approaches infinity. However, in
practice, even  a system is stable, it may be totally useless,
because it possesses unsatisfactory transient performance. So, one
should  be interested in not only classical Lyapunov stability, but
also finite-time transient performance. To this end, finite-time
stability was proposed in 1950s \cite{Dorato,Kamenkov}. Finite-time
stability is
  different from classical stability in two aspects: First, the
concerned system operation is confined to a prescribed  finite-time
interval instead of an infinite-time horizon. Second, the state
trajectory lies within a specific bound over the given finite
interval of time.  Recently, finite-time stability has become a
popular research topic due to its practical sense. Many nice results
have been obtained, and we refer the reader to \cite{Amato,17} for
linear deterministic systems,
 stochastic It\^o    systems
\cite{10,quanxin,yin2011,zhang-an,newyan1,newyan2}, switching
systems \cite{newguopei,17,8,yanpark}. It is worth mentioning that
short-time stability is also referred to finite-time stability as
discussed in \cite{10,Lebedev,quanxin,yin2011}. The reference
\cite{zhang-an} extended  finite-time stability and stabilization of
\cite{Amato} to linear time-invariant It\^o systems.    A
mode-dependent parameter approach was proposed to give a sufficient
condition for finite-time stability and stabilization for It\^o
stochastic systems with Markovian switching \cite{8}, and the same
method was also used to deal with finite-time guaranteed cost
control of It\^o stochastic Markovian jump systems\cite{yanpark}. It
can be found that, up to now, there are few  results on finite-time
stability and stabilization of discrete  time-varying stochastic
systems with multiplicative noise.

As said by J.P. LaSalle  \cite{laSallebook}, ``Today there is more
and more reason for studying difference equations systematically.
They are in their own right important mathematical models",
``Moreover, their study provides  a good introduction  to the study
theory of differential equations, difference-differential equations,
and functional differential equations". It is expected that, along
the development of computer techniques, the study  on  discrete-time
systems will become  more and more important, and attract a lot of
researchers' attention. In \cite{newNi1,newNi2,newNi3,newNi4},
discrete-time mean-field linear-quadratic optimal control problems
have been systematically researched.  In
\cite{dragan,Bouhtouri_1999,book}, the $H_\infty$ control of linear
discrete-time stochastic systems were studied, especially, necessary
and sufficient conditions for mean square stability were presented.
In \cite{5}, robust stability and stabilization  for  a class of
linear  discrete-time time-varying stochastic systems   with
Markovian jump were investigated based on a small-gain theorem. The
reference \cite{6} was about finite-time stability of discrete
time-varying systems with randomly occurring nonlinearity and
missing measurements. Finite-time stochastic stability and
stabilisation with partly unknown transition probabilities for
linear discrete-time Markovian jump systems was considered in
\cite{11}. By choosing Lyapunov-Krasovskii-like functionals,
sufficient conditions were given in \cite{14,13} for finite-time
stability of linear deterministic systems with time-varying delay
based on LMIs. Most results on finite-time stability of stochastic
systems are sufficient but not necessary conditions, which are
derived by Lyapunov function or Lyapunov-Krasovskii-like  functional
method.

This  paper will study the finite-time stability of the following
linear discrete time-varying stochastic system with multiplicative
noise
$$
x_{k+1}=A_kx_{k}+C_kx_{k}w_k,
$$
and the finite-time stabilization of the following control system
$$
x_{k+1}=A_kx_{k}+B_ku_{k}+(C_kx_{k}+D_ku_k)w_k.
$$
Note that in order to study the detectability and observability of the following linear discrete time-varying system
\begin{equation}\label{bvhb}
\begin{cases}
 \  x_{k+1}=A_kx_{k}+C_kx_{k}w_k, \\
 \ y(k)=H_kx(k),
   \end{cases}
\end{equation}
we introduce the  state transition matrix of (\ref{bvhb}) in mean
square sense  \cite{IJC2016}, while this paper can be viewed as the
continuation of \cite{IJC2016}. The contributions of this paper are
as follows:

(i)  We develop a state transition matrix approach to present some necessary and sufficient conditions for finite-time
stability of linear discrete stochastic systems. Specifically,  the following two identities
$$
\phi_{l,k}'(I_{2^{l-k}}\otimes
R_l)\phi_{l,k}=R_k^{\frac{1}{2}}\bar{\phi}_{l,k}'\bar{\phi}_{l,k}R_k^{\frac{1}{2}}
$$
and
$$
\psi_{l,k}'(I_{2^{l-k}}\otimes
R_l)\psi_{l,k}=R_k^{\frac{1}{2}}\bar{\psi}_{l,k}'\bar{\psi}_{l,k}R_k^{\frac{1}{2}}
$$
proved in this paper   are   important, which have potential
applications to the study of piecewise finite-time stability and
other control issues.

(ii)  In order to  further study  finite-time  stabilization and
obtain easily testing criteria, we apply the Lyapunov function
method to present some sufficient conditions for finite-time
stability and stabilization based on matrix inequalities and   LMIs.

The paper is organized as follows: In Section 2, we define
finite-time stability and stabilization for linear  discrete
time-varying stochastic systems. Several useful lemmas are
presented.   In Section 3,  finite-time stability is studied based
on state transition matrix approach, where some  necessary and
sufficient criteria are obtained for finite-time stability.   In
Section 4, we make use of Lyapunov function method to investigate
finite-time stability and stabilization,  and several sufficient
conditions for  finite-time stability and stabilization are given
based on matrix inequalities and LMIs. In Section 5, two examples
are constructed to show the effectiveness of our obtained results.

For convenience, the notations adopted in this paper are as follows.

$M'$: the transpose of the matrix $M$ or vector $M$; $M>0$ ($M<0$):
the matrix $M$ is a positive definite (negative definite) symmetric
matrix; $I_n$: $n\times n$ identity matrix; ${\mathcal R}^n$: the
$n$-dimensional real Euclidean vector space; ${\mathcal R}^{n\times
m}$: the space of all $n\times m$ matrices with entries in
${\mathcal R}$; $A\otimes B$: the Kronecker product of two matrices
$A$ and $B$. $\lambda_{min}(A) (\lambda_{max}(A)$): the minimum
(maximum) eigenvalue of a real symmetric  matrix $A$;
 ${\mathcal N}:=\{0,1,2,\cdots,\}$;  $N_+:=\{1,2,\cdots,\}$;  $N_T:=\{0,1,\cdots, T\}$ where   $T\in N_+$.

\section{Preliminaries\label{sec:PR}}
\hspace{0.13in}Consider the following discrete-time time-varying
stochastic difference system  described by
\begin{equation}
\left\{\begin{array}{cccc}
\label{system1}
 \  x_{k+1}=A_kx_{k}+C_kx_{k}w_k, \\
 \ x_0\in {\mathcal R}^n, k\in {N_{T-1}}.
   \end{array}
   \right.
\end{equation}
where $x_{k}\in \mathcal{R}^n$ is the $n$-dimensional state vector.
$\{w_k\}_{k\in N_{T-1}}$ is a sequence of  one-dimensional independent
white noise processes  defined on the complete  filtered probability
space $(\Omega, {\mathcal F}, {\{{\mathcal F}_k\}}_{k\in N_{T}},
{\mathcal\mathcal P})$, where ${\mathcal
F}_k=\sigma(w_0,w_1,\cdots,w_{k-1})$, ${\mathcal
F}_0=\{\phi,\Omega\}$.   Assume that $E[w_k]=0$,
$E[w_kw_j]=\delta_{kj}$, where $E$ stands for the mathematical expectation operator,  and
$\delta_{kj}$ is a Kronecker function defined by $\delta_{kj}=0$ for
$k\ne j$ while $\delta_{kj}=1$ for $k=j$. Without loss of
generality,   $x_0$ is assumed to be a deterministic vector. $A_k$ and $C_k$ are $n\times n$ time-varying matrices with respect to $k$.

{\bf Definition 1}\label{def1}
Given  a positive integer $T$, two positive scalars $0<c_1\le c_2$,
and a finite positive difinite symmetric matrix  sequence
$\{R_k>0\}_{k\in N_T}$. The system (\ref{system1}) is said to be
finite-time stable with respect to $(c_1, c_2, T, \{R_k\}_{k\in
N_T})$, if
\begin{equation}
\|x_0\|^2_{R_0}\leq c_1\Rightarrow E\|x_{k}\|^2_{R_k}<c_2, \ \ \
\forall k \in N_{T},
\end{equation}
where  $\|x\|^2_R:=x'Rx$.

This paper will also study the feedback stabilization of the following linear discrete  time-varying control system
\begin{equation}
\left\{\begin{array}{cccc}
\label{system3}
 \  x_{k+1}=A_kx_{k}+B_ku_{k}+(C_kx_{k}+D_ku_k)w_k, \\
 \ x_0\in {\mathcal R}^n, k\in {N_{T-1}},
   \end{array}
   \right.
\end{equation}
where $u_k\in \mathcal{R}^m$ is the $m$-dimensional control input.

{\bf Definition 2}\label{def2}
System (\ref{system3}) is said to be finite-time stabilizable with
respect to $(c_1, c_2, T, \{R_k\}_{k\in {N_{T}}})$, if there exists
a linear  state feedback control law $u_k=K_kx_k$, $k\in N_{T-1}$,
such that  the resulting closed-loop system
\begin{equation}
\left\{\begin{array}{cccc}
\label{system4}
 \  x_{k+1}=(A_k+B_kK_k)x_k+(C_k+D_kK_k)x_{k}w_k, \\
 \ x_0\in {\mathcal R}^n, \ k\in {N_{T-1}}
   \end{array}
   \right.
\end{equation}
is finite-time stable with respect to $(c_1, c_2, T, \{R_k\}_{k\in N_T})$.
In order to investigate the finite-time stability of system  (\ref{system1}), we need to introduce some useful lemmas.

{\bf Lemma 1} \cite{IJC2016}\label{lem1}\
For system (\ref{system1}), we have
\begin{itemize}
\item[(i)] $E\|x_l\|^2=E\|\phi_{l,k}x_{k}\|^2$ for $ l\ge k$, where
$\phi_{j,j}=I_n$ for $j\in {\mathcal N}$, and
\begin{equation}\label{state2}
\phi_{l,k}=\left[
\begin{array}{cc}
(I_{2^{l-k-1}}\otimes A_{l-1})\phi_{l-1,k}\\
(I_{2^{l-k-1}}\otimes C_{l-1})\phi_{l-1,k}
\end{array}
\right],\quad l>k.
\end{equation}
\item[(ii)] $E\|x_l\|^2=E\|\psi_{l,k}x_{k}\|^2$ for $ l\ge k$, where
$\psi_{j,j}=I_n$ for $j\in {\mathcal N}$,  and
\begin{equation}\label{state1}
\psi_{l,k}=\left[
\begin{array}{cc}
\psi_{l,k+1}A_k\\
\psi_{l,k+1}C_k
\end{array}
\right],\quad l>k.
\end{equation}
\item[(iii)] $x_{k}\in l^2_{{\mathcal
F}_{k}}$ if $A_i$ and $C_i$ are uniformly bounded for $i\in {\mathcal N}$.
\end{itemize}

{\bf Remark 1}
The matrices  $\phi_{l,k}$ and $\psi_{l,k}$ defined in Lemma 1 can be
viewed  as the state transition matrices  in the mean square sense.
In \cite{IJC2016}, $\phi_{l,k}$ and $\psi_{l,k}$ were introduced to present  detectability conditions  of linear discrete time-varying
stochastic systems.

From Lemma 1,  the  state transition matrix  is not unique, since that
$\psi_{l,k}=\phi_{l,k}$ is not necessarily true; see the following example.

{\bf Example 1}
It is easy to compute that
$$
\psi_{l,l-1}=\phi_{l,l-1}
=\left[\begin{array}{ccc}
A_{l-1}\\
C_{l-1}
\end{array}
\right],
$$
$$
\psi_{l,l-2}=
\left[\begin{array}{c}
A_{l-1}A_{l-2}\\
C_{l-1}A_{l-2}\\
A_{l-1}C_{l-2}\\
C_{l-1}C_{l-2}
\end{array}
\right], \
\phi_{l,l-2}=
\left[\begin{array}{c}
A_{l-1}A_{l-2}\\
A_{l-1}C_{l-2}\\
C_{l-1}A_{l-2}\\
C_{l-1}C_{l-2}
\end{array}
\right].
$$
Hence, $\phi_{l,l-2}\ne\psi_{l,l-2}$ for $A_{l-1}\ne C_{l-1}$, but
we always have $\phi'_{l,k}\phi_{l,k}=\psi'_{l,k}\psi_{l,k}$ for
$l\ge k$.

{\bf Lemma 2}(Schur's complement)\label{schur's}
\ \ For a real symmetric matrix
$S=\left[\begin{array}{ccc}S_{11}&S_{12}\\S_{12}^T&S_{22}\end{array}\right]$,
 the following three
conditions are equivalent:
\begin{itemize}
\item [(i)]$S<0$;
\item [(ii)]$S_{11}<0, S_{22}-S^T_{12}S^{-1}_{11}S_{12}<0$;
\item [(iii)] $S_{22}<0, S_{11}-S_{12}S^{-1}_{22}S^T_{12}<0$.
\end{itemize}

{\bf Lemma 3}\cite{Horn}\label{lem5}
For matrices $A$, $B$, $C$ and $D$ of suitable dimensions, we
have
\begin{equation}
(A\otimes B)(C\otimes D)=(AC)\otimes(BD).
\end{equation}
\section{State Transition Matrix-based Approach for Finite-time Stability \label{sec:PARETOD}}

\hspace{0.13in}  In this section, we mainly use the state transition
matrix  approach developed in \cite{book,IJC2016} to study the
finite-time stability and stabilization of the system
(\ref{system1}).

Set  $\bar{x}_{k}=R_{k}^{\frac{1}{2}}x_{k}$, then
$E[x_k'R_kx_k]=E[\bar{x}_k'\bar{x}_k]$, and  system (\ref{system1})
is equivalent to  the following system:
\begin{eqnarray}\label{system5}
\left\{\begin{array}{cccc}
 \  \bar{x}_{k+1}=R_{k+1}^{\frac{1}{2}}A_kR_{k}^{-\frac{1}{2}}\bar{x}_{k}+R_{k+1}^{\frac{1}{2}}C_kR_{k}^{-\frac{1}{2}}\bar{x}_{k}w_k, \\
 \ \ \ \ \ \ \ \ \  =\bar{A}_k\bar{x}_{k}+\bar{C}_k \bar{x}_{k}w_k,\\
 \ \bar{x}_0=R_0^{\frac{1}{2}}x_0\in {\mathcal R}^n, k\in N_{T-1},
   \end{array}
   \right.
\end{eqnarray}
where $\bar{A}_k=R_{k+1}^{\frac{1}{2}}A_kR_{k}^{-\frac{1}{2}}$,
$\bar{C}_k=R_{k+1}^{\frac{1}{2}}C_kR_{k}^{-\frac{1}{2}}$. Lemma 1 yields directly  the following lemma:

{\bf Lemma 4}\label{lem3.1}
The state transition matrix $\bar{\phi}_{l,k}$ for the  system
(\ref{system5}) is in the form of
\begin{eqnarray}\label{newstate1}
\bar{\phi}_{l,k}
&=&\left[
\begin{array}{cc}
(I_{2^{l-k-1}}\otimes \bar{A}_{l-1})\bar{\phi}_{l-1,k}\\
(I_{2^{l-k-1}}\otimes \bar{C}_{l-1})\bar{\phi}_{l-1,k}
\end{array}
\right], \ l>k, \ \bar{\phi}_{k,k}=I.
\end{eqnarray}
Moreover,  $E\|\bar{x}_l\|^2=E\|\bar{\phi}_{l,k}\bar{x}_{k}\|^2$.
Note that $E\|x_k\|^2_{R_k}=E\|\bar{x}_k\|^2$, we can obtain the
following result directly.

{\bf Proposition 1}\label{PRO3.1}
The system (\ref{system1}) is finite-time stable with respect to
$(c_1, c_2, T, \{R_k\}_{k\in N_T})$ if and only if (iff) the system
(\ref{system5}) is finite-time stable with respect to $(c_1, c_2, T,
\{R_k=I_n\}_{k\in N_T})$.
The next lemma establishes the relationship between $\phi_{k,0}$ and
$\bar{\phi}_{k,0}$.

{\bf Lemma 5}\label{lem3.2} The state transition matrix $\phi_{k,0}$
of system  (\ref{system1})  and the state transition matrix
$\bar{\phi}_{k,0}$ of system (\ref{system5}) have the following
relation:
\begin{equation} \label{ass1}
\phi_{k,0}'(I_{2^k}\otimes
R_k)\phi_{k,0}=R_0^{\frac{1}{2}}\bar{\phi}_{k,0}'\bar{\phi}_{k,0}R_0^{\frac{1}{2}}.
\end{equation}

{\textbf {Proof.}} Lemma 5  can be shown by induction.
Because
$$
\phi_{k,0}'(I_{2^k}\otimes R_k)\phi_{k,0}=\phi_{k,0}'(I_{2^k}\otimes
R^{1/2}_k)(I_{2^k}\otimes R^{1/2}_k)\phi_{k,0},
$$
we only need to prove the following identity:
\begin{equation}\label{efra}
 (I_{2^k}\otimes
R_k^{\frac{1}{2}})\phi_{k,0}=\bar{\phi}_{k,0}R_0^{\frac{1}{2}}.
\end{equation}
For $k=0$, in view of  $\phi_{0,0}=\bar{\phi}_{0,0}=I_n$, we have
$$
(I_1\otimes
R_0^{\frac{1}{2}})\phi_{0,0}=R_0^{\frac{1}{2}}=\bar{\phi}_{0,0}R_0^{\frac{1}{2}}.
$$
Hence, (\ref{efra}) holds for $k=0$. Assume that for $k=j-1$,
(\ref{efra}) holds, i.e., $(I_{2^{j-1}}\otimes
R_{j-1}^{\frac{1}{2}})\phi_{{j-1},0}=\bar{\phi}_{{j-1},0}R_0^{\frac{1}{2}}$,
then   we  shall prove $(I_{2^j}\otimes
R_j^{\frac{1}{2}})\phi_{j,0}=\bar{\phi}_{j,0}R_0^{\frac{1}{2}}$. It
can be seen that
\begin{eqnarray*}
(I_{2^j}\otimes R_j^{\frac{1}{2}})\phi_{j,0}&=&
\left[\begin{array}{ccc} I_{2^{j-1}}\otimes R_j^{\frac{1}{2}} & \ \
\ \  0\\
0 &  I_{2^{j-1}}\otimes R_j^{\frac{1}{2}}
\end{array}\right]\left[
\begin{array}{cc}
(I_{2^{j-1}}\otimes A_{j-1})\phi_{j-1,0}\\
(I_{2^{j-1}}\otimes C_{j-1})\phi_{j-1,0}
\end{array}
\right]\\
&=&\left[\begin{array}{cc}
(I_{2^{j-1}}\otimes (R_j^{\frac{1}{2}}A_{j-1}))\phi_{j-1,0}\\
(I_{2^{j-1}}\otimes (R_j^{\frac{1}{2}}C_{j-1}))\phi_{j-1,0}
\end{array}
\right].
\end{eqnarray*}
By Lemma 3 and (\ref{state2}), we have {\small
\begin{eqnarray}
&&(I_{2^j}\otimes
R_j^{\frac{1}{2}})\phi_{j,0}\nonumber\\
&=&\left[\begin{array}{cc}
(I_{2^{j-1}}\otimes (R_j^{\frac{1}{2}}A_{j-1}))\phi_{j-1,0}\\
(I_{2^{j-1}}\otimes (R_j^{\frac{1}{2}}C_{j-1}))\phi_{j-1,0}
\end{array}
\right]\nonumber\\
&=&\left[\begin{array}{cc}
(I_{2^{j-1}}\otimes (R_j^{\frac{1}{2}}A_{j-1}))(I_{2^{j-1}}\otimes R_{j-1}^{-\frac{1}{2}})(I_{2^{j-1}}\otimes R_{j-1}^{\frac{1}{2}})\phi_{j-1,0}\\
(I_{2^{j-1}}\otimes (R_j^{\frac{1}{2}}C_{j-1}))(I_{2^{j-1}}\otimes
R_{j-1}^{-\frac{1}{2}})(I_{2^{j-1}}\otimes
R_{j-1}^{\frac{1}{2}})\phi_{j-1,0}
\end{array}
\right]\nonumber\\
&=&\left[\begin{array}{cc} (I_{2^{j-1}}\otimes
(R_j^{\frac{1}{2}}A_{j-1}R_{j-1}^{-\frac{1}{2}}))(I_{2^{j-1}}\otimes
R_{j-1}^{\frac{1}{2}})\phi_{j-1,0}\\(I_{2^{j-1}}\otimes
(R_j^{\frac{1}{2}}C_{j-1}R_{j-1}^{-\frac{1}{2}}))(I_{2^{j-1}}\otimes
R_{j-1}^{\frac{1}{2}})\phi_{j-1,0}
\end{array}
\right]\nonumber\\
&=&\bar{\phi}_{j,0}R_0^{\frac{1}{2}}.
\end{eqnarray}}
The proof of this lemma is completed. $\square$

Repeating the same procedure as in Lemma 5, the following
more general relation still holds.

{\bf Lemma 6}\label{lem3.2ff} For any $l\ge k$, $l, k\in {\mathcal N}$,
we have the following identity:
\begin{equation} \label{ass1ccaa}
\phi_{l,k}'(I_{2^{l-k}}\otimes
R_l)\phi_{l,k}=R_k^{\frac{1}{2}}\bar{\phi}_{l,k}'\bar{\phi}_{l,k}R_k^{\frac{1}{2}}.
\end{equation}
Corresponding to the second form of the state transition matrix $\psi_{l,k}$ defined in   (\ref{state1}),
we denote $\bar{\psi}_{l,k}$ as another  state transition matrix of system  (\ref{system1}) from the state $x_k$
to the state $x_l$ with $l\ge k$, which is defined by
$$
\bar{\psi}_{l,k}=\left[
\begin{array}{cc}
\bar{\psi}_{l,k+1}\bar{A}_k\\
\bar{\psi}_{l,k+1}\bar{C}_k
\end{array}
\right], \ \bar{\psi}_{l,l}=I.
$$
Now, a more general relation corresponding to (\ref{ass1ccaa}) can
be given based on $\bar{\psi}_{l,k}$.

{\bf Lemma 7}\label{lem3.2rr} The state transition matrix $\psi_{l,k}$
of system  (\ref{system1})  and the state transition matrix
$\bar{\psi}_{l,k}$ of system (\ref{system5}) have the following
relation:
\begin{equation} \label{ass1cc}
\psi_{l,k}'(I_{2^{l-k}}\otimes
R_l)\psi_{l,k}=R_k^{\frac{1}{2}}\bar{\psi}_{l,k}'\bar{\psi}_{l,k}R_k^{\frac{1}{2}}.
\end{equation}

{\textbf {Proof.}} We still prove this lemma by induction. For
$k=l$,
$$
\psi_{l,k}'(I_{2^{l-k}}\otimes
R_l)\psi_{l,k}=\psi_{l,l}'(I_{2^{0}}\otimes
R_l)\psi_{l,l}=R_{l},\\
$$
while
$$
R_k^{\frac{1}{2}}\bar{\psi}_{l,k}'\bar{\psi}_{l,k}R_k^{\frac{1}{2}}
=R_l^{\frac{1}{2}}\bar{\psi}_{l,l}'\bar{\psi}_{l,l}R_l^{\frac{1}{2}}
=R_l.
$$
So
in the case $k=l$, (\ref{ass1cc}) holds. Assume for $k=j$, (\ref{ass1cc}) holds, i.e.,
\begin{equation} \label{agfg1cc}
\psi_{l,j}'(I_{2^{l-j}}\otimes
R_l)\psi_{l,j}=R_j^{\frac{1}{2}}\bar{\psi}_{l,j}'\bar{\psi}_{l,j}R_j^{\frac{1}{2}}.
\end{equation}
Then we only need to show that (\ref{ass1cc}) holds for $k=j-1$,
i.e.,
\begin{eqnarray} \label{agfdfg1cc}
\psi_{l,j-1}'(I_{2^{l-j+1}}\otimes
R_l)\psi_{l,j-1}=R_{j-1}^{\frac{1}{2}}\bar{\psi}_{l,j-1}'\bar{\psi}_{l,j-1}R_{j-1}^{\frac{1}{2}}.
\end{eqnarray}
By induction assumption (\ref{agfg1cc}), the right hand side of  (\ref{agfdfg1cc}) can be computed as
\begin{eqnarray*}
&&R_{j-1}^{1/2}\bar{\psi}_{l,j-1}'\bar{\psi}_{l,j-1}R_{j-1}^{1/2}\\
&=&R_{j-1}^{1/2} \left[
\begin{array}{ccc}
R_{j-1}^{-1/2}A'_{j-1}R_j^{1/2}\bar{\psi}_{l,j}' &
R_{j-1}^{-1/2}C'_{j-1}R_j^{1/2}\bar{\psi}_{l,j}'
\end{array}
\right]\\
&&\cdot \left[
\begin{array}{ccc}
 \bar{\psi}_{l,j}R_j^{1/2}A_{j-1}R_{j-1}^{1/2}\\
 \bar{\psi}_{l,j}R_j^{1/2}C_{j-1}R_{j-1}^{1/2}
\end{array}
\right]R_{j-1}^{1/2}\\
&=&
\left[
\begin{array}{ccc}
A'_{j-1}R_j^{1/2}\bar{\psi}_{l,j}' & C'_{j-1}R_j^{1/2}\bar{\psi}_{l,j}'
\end{array}
\right]
\left[
\begin{array}{ccc}
 \bar{\psi}_{l,j}R_j^{1/2}A_{j-1}\\
 \bar{\psi}_{l,j}R_j^{1/2}C_{j-1}
\end{array}
\right]\\
&=&\left[
\begin{array}{ccc}
A'_{j-1} & C'_{j-1}
\end{array}
\right]\left[
\begin{array}{ccc}
R_j^{1/2}\bar{\psi}_{l,j}\bar{\psi}_{l,j}R_j^{1/2} & 0\\
0 & R_j^{1/2}\bar{\psi}_{l,j}\bar{\psi}_{l,j}R_j^{1/2}
\end{array}
\right] \left[
\begin{array}{ccc}
A_{j-1}\\
C_{j-1}
\end{array}
\right]\\
&=&\psi_{l,j-1}'(I_{2^{l-j+1}}\otimes
R_l)\psi_{l,j-1}.
\end{eqnarray*}
Hence, (\ref{agfdfg1cc}) holds, this lemma is proved. $\square$

Lemmas 6-7 have potential important
applications to   piecewise finite-time stability and mean square
stability.

{\bf Theorem 1}\label{th1} The system (\ref{system1})
is finite-time stable with respect to $(c_1, c_2, T, \{R_k\}_{k\in
N_T})$  iff  the following inequalities are satisfied:
\begin{equation} \label{1as}
\phi_{k,0}'(I_{2^k}\otimes R_k)\phi_{k,0}<\frac{c_2}{c_1}R_0,\ \ \
\forall k\in N_T.
\end{equation}

{\textbf {Proof.}}  We first prove the sufficiency of Theorem
1. If
\begin{equation}\label{dtfhj}
x_0'R_0x_0 \leq c_1,
\end{equation}
then, by Lemma 5, we have
\begin{eqnarray}\label{eq fgyh}
E\|x_{k}\|^2_{R_k}
&=&E\|{\bar x}_{k}\|^2={\bar x}'_0\bar{\phi}_{k,0}'\bar{\phi}_{k,0}{\bar x}_0
=x_0'R^{1/2}_0\bar{\phi}_{k,0}'\bar{\phi}_{k,0}R^{1/2}_0{x}_0\nonumber\\
&=& x_0'\phi_{k,0}'(I_{2^k}\otimes R_k)\phi_{k,0}x_0
=\|x_0\|^2_{\phi_{k,0}'(I_{2^k}\otimes R_k)\phi_{k,0}}.
\end{eqnarray}
If $x_0=0$, then
\begin{equation}\label{hjhj}
E\|x_{k}\|^2_{R_k}\equiv 0<c_2, \forall k\in N_T.
\end{equation}
If $x_0\ne 0$, then, by (\ref{dtfhj}) and (\ref{eq fgyh}), it
follows that
\begin{equation}\label{hddfsdw}
E\|x_{k}\|^2_{R_k}<\|x_0\|^2_{\frac{c_2}{c_1}R_0}\le c_2, \ k\in
{N_T}.
\end{equation}
In general, for any $x_0\in {R}^n$, we can always conclude
$E\|x_{k}\|^2_{R_k}<c_2$ from $x_0'R_0x_0 \leq c_1$, i.e., the
system (\ref{system1}) is finite-time stable with respect to $(c_1,
c_2, T, \{R_k\}_{k\in N_T})$.

Below,  we prove the necessity part by contradiction. If  the system
(\ref{system1}) is finite-time stable with respect to $(c_1, c_2, T,
\{R_k\}_{k\in N_T})$,  while there exists  $\bar{k}\in {N_T}$ such that
\begin{equation} \label{1assds}
\phi_{\bar{k},0}'(I_{2^{\bar{k}}}\otimes R_{\bar{k}})\phi_{\bar{k},0}\ge
\frac{c_2}{c_1}R_0,
\end{equation}
then by taking an $x_0$ satisfying $x_0'R_0x_0=c_1$, it yields that
\begin{eqnarray}\label{eq cgcx}
E\|{x}_{\bar k}\|^2_{R_{\bar
k}}&=&E\|{x_0}\|^2_{\phi_{\bar{k},0}'(I_{2^{\bar{k}}}\otimes
R_{\bar{k}})\phi_{\bar{k},0}}\geq {x_0}'\frac{c_2}{c_1}R_0{x_0}\geq c_2.
\end{eqnarray}
The inequality (\ref{eq cgcx}) contradicts the finite-time
stability, which requires
$$
E\|{x}_{\bar k}\|^2_{R_{\bar k}}<c_2.
$$
The proof is completed. $\square$

{\bf Corollary 1}\label{thnew1}
The system (\ref{system1}) is finite-time stable with respect to
 $(c_1, c_2, T, \{R_k\}_{k\in {N_T}})$ if the following holds for $k\in
 N_T$.
\begin{equation}\label{21}
\phi_{k,0}'\phi_{k,0}<\frac{c_2R_0}{c_1\max_{k\in
N_T}\lambda_{\max}(R_k)},
\end{equation}
where $\lambda_{max}(R_k)$ denotes  the  maximum eigenvalue of the
 matrix  $R_k$.

{\textbf {Proof.}} Note that
$$
\phi_{k,0}'(I_{2^k}\otimes R_k)\phi_{k,0}\leq \max_{k\in
N_T}\lambda_{\max}(R_k) {\phi_{k,0}}'{\phi_{k,0}}.
$$
Hence, by Theorem 1,  the condition (\ref{21}) is a
sufficient condition for the finite-time stability of the system
(\ref{system1}). $\square$

Because $R_k>0$, $k\in {N_T}$,  there is a nonsingular matrix
sequence $\{L_k\}_{k\in N_T}$, such that $R_k=L_k'L_k$. By Lemma
3, we obtain ${I_{2^k}\otimes R_k}=(I_{2^k}\otimes
L_k')(I_{2^k}\otimes L_k)$. By Theorem 1, the solvability of
inequality (\ref{1as}) is equivalent to that  the system
(\ref{system1}) is finite-time stable with respect to $(c_1, c_2, T,
\{R_k\}_{k\in {N_T}})$. In addition, for $\forall k\in N_T$, we have
\begin{eqnarray}\label{lass}
&&\phi_{k,0}'(I_{2^k}\otimes R_k)\phi_{k,0}< \frac{c_2}{c_1}R_0\nonumber\\
\Longleftrightarrow&&\phi_{k,0}'(I_{2^k}\otimes L_k')(I_{2^k}\otimes
L_k)\phi_{k,0}-\frac{c_2}{c_1}R_0<0\nonumber\\
\Longleftrightarrow &&R_0^{-\frac{1}{2}}\phi_{k,0}'(I_{2^k}\otimes
L_k')(I_{2^k}\otimes L_k)\phi_{k,0}R_0^{-\frac{1}{2}} -\frac{c_2}{c_1}I_n<0.
\end{eqnarray}
Using Lemma 2 twice,  the inequality (\ref{lass}) can be
equivalently written as
\begin{eqnarray}\label{fhgdsj}
&&\left[\begin{array}{ccc}
-\frac{c_2}{c_1}I_n&R_0^{-\frac{1}{2}}\phi_{k,0}'(I_{2^k}\otimes
L_k')\\
(I_{2^k}\otimes
L_k)\phi_{k,0}R_0^{-\frac{1}{2}}&-I_{2^kn}\end{array}\right]<0
\nonumber\\
\Longleftrightarrow && (I_{2^k}\otimes
L_k)\phi_{k,0}R_0^{-\frac{1}{2}}R_0^{-\frac{1}{2}}\phi_{k,0}'(I_{2^k}\otimes
L_k')-\frac{c_2}{c_1}I_{2^kn}<0.
\end{eqnarray}
Let $P_k=\phi_{k,0}R_0^{-1}\phi_{k,0}'$, then (\ref{fhgdsj}) leads
to
\begin{eqnarray}\label{erwwy}
(I_{2^k}\otimes L_k)P_k(I_{2^k}\otimes
L_k')-\frac{c_2}{c_1}I_{2^kn}<0.
\end{eqnarray}
Pre-multiplying $(I_{2^k}\otimes L_k)^{-1}$ and post-multiplying
$[(I_{2^k}\otimes L_k)^{-1}]^{'}$ on  both sides of (\ref{erwwy}),
it follows that
\begin{eqnarray}
P_k<\frac{c_2}{c_1}(I_{2^k}\otimes R_k^{-1}).
\end{eqnarray}
Recall the following properties of $\phi_{\cdot,\cdot}$ in Lemma
1:
\begin{eqnarray}
\phi_{0,0}=I,\ \
\phi_{k+1,0}=\left[
\begin{array}{cc}
(I_{2^{k}}\otimes A_k)\phi_{k,0}\\
(I_{2^{k}}\otimes C_k)\phi_{k,0}
\end{array}
\right],
\end{eqnarray}
if we substitute  $\phi_{k+1,0}$ into
$P_{k+1}=\phi_{k+1,0}R_0\phi_{k+1,0}'$, then {\small
\begin{eqnarray}
&&P_{k+1}=\left[
\begin{array}{cc}
(I_{2^{k}}\otimes A_k)\phi_{k,0}\\
(I_{2^{k}}\otimes C_k)\phi_{k,0}
\end{array}
\right]R_0^{-1} \left[
\begin{array}{cc}
\phi_{k,0}'(I_{2^{k}}\otimes A_k')& \phi_{k,0}'(I_{2^{k}}\otimes
C_k')
\end{array}
\right]\nonumber\\
&&=\left[
\begin{array}{ccc}
(I_{2^{k}}\otimes A_k)P_k(I_{2^{k}}\otimes A_k')& (I_{2^{k}}\otimes
A_k)P_k(I_{2^{k}}\otimes C_k')
\\(I_{2^{k}}\otimes
C_k)P_k(I_{2^{k}}\otimes A_k')&(I_{2^{k}}\otimes
C_k)P_k(I_{2^{k}}\otimes C_k')
\end{array}
\right].
\end{eqnarray}
} Summarize  the previous discussion, we are in a position to obtain
the following theorem, which is  equivalent to  Theorem 1.

{\bf  Theorem 2}\label{th2}
The system (\ref{system1}) is finite-time stable with respect to
$(c_1, c_2, T, \{R_k\}_{k\in N_T})$ iff there exists a symmetric
matrix  sequence $\{P_k\}_{k\in\{0,1,\cdots,T\}}$  solving the
following constrained difference equation: {\small
\begin{eqnarray}
\left\{\begin{array}{l}
P_{0}=R_{0}^{-1}, \\
P_{k+1}=\left[
\begin{array}{cc}
(I_{2^{k}}\otimes A_k)P_k(I_{2^{k}}\otimes A_k')& (I_{2^{k}}\otimes
A_k)P_k(I_{2^{k}}\otimes C_k')
\\(I_{2^{k}}\otimes
C_k)P_k(I_{2^{k}}\otimes A_k')&(I_{2^{k}}\otimes
C_k)P_k(I_{2^{k}}\otimes C_k') \nonumber
\end{array}
\right],\\
P_k<\frac{c_2}{c_1}(I_{2^k}\otimes R_k^{-1}),\label{eq vfhgg}\ k\in
N_T.
\end{array}
\right.
\end{eqnarray}}
If  system (\ref{system1})  is finite-time stable with respect to
$(c_1, c_2, T, \{R_k\}_{k\in N_T})$, we shall show that  the
following  perturbed system
\begin{equation}\label{system2}
\left\{\begin{array}{cccc}
 \  x_{k+1}=(A_k+\varepsilon I)x_{k}+(C_k+\varepsilon I)x_{k}w_k, \\
 \ x_0\in {\mathcal R}^n, k\in {N}_{T-1},
   \end{array}
   \right.
\end{equation}
is also finite-time stable with respect to $(c_1$, $c_2$, $T$,
$\{R_k\}_{k\in N_T})$  for sufficiently small real number
$\varepsilon$.

{\bf Theorem 3}\label{the3.3}
If  system (\ref{system1})  is finite-time stable with respect to
$(c_1, c_2, T, \{R_k\}_{k\in N_T})$, then so is the system
(\ref{system2}).

{\textbf {Proof.}} From Theorem 2, system (\ref{system2})
 is finite-time stable with respect to $(c_1, c_2, T,
\{R_k\}_{k\in N_T})$  iff    the constrained difference equation is
admissible.{\small
\begin{eqnarray}\label{21221}
\begin{array}{l}
P^{\varepsilon}_{0}=R_{0}^{-1}, \\
P^{\varepsilon}_{k+1}=\left[
\begin{array}{cccccc}
[I_{2^{k}}\otimes (A_k+\varepsilon
I)]P^{\varepsilon}_k[I_{2^{k}}\otimes (A_k+\varepsilon I)']
&[I_{2^{k}}\otimes (A_k+\varepsilon I)]P^{\varepsilon}_k[I_{2^{k}}\otimes
(C_k+\varepsilon I)']
\\
*
&
 [I_{2^{k}}\otimes (C_k+\varepsilon
I)]P^{\varepsilon}_k[I_{2^{k}}\otimes (C_k+\varepsilon I)']
\end{array}\right],\\
P^\varepsilon_k<\frac{c_2}{c_1}(I_{2^k}\otimes R_k^{-1}), k\in N_T\\
\end{array}
\end{eqnarray}}
From (\ref{21221}), we can see that $P^\varepsilon_k\to P_k$ as
$\varepsilon\to 0$, where $P_k$ is defined by (\ref{eq vfhgg}). If
system (\ref{system1})  is finite-time stable with respect to $(c_1,
c_2, T, \{R_k\}_{k\in N_T})$, then,  by Theorem 2,
$P_k<\frac{c_2}{c_1}(I_{2^k}\otimes R_k^{-1})$, $k\in N_T$, which
implies $P^\varepsilon_k<\frac{c_2}{c_1}(I_{2^k}\otimes R_k^{-1})$
for sufficiently small $\varepsilon>0$. This theorem is proved.
$\square$

\section{Lyapunov Function-based Approach for Finite-time Stability and Stabilization \label{sec:LMI}}

\hspace{0.13in} The state transition matrix-based approach presents
necessary and sufficient conditions for finite-time stability of
system (\ref{system1}), which is elegant in theory.   However, there
is some difficulty  in applying the state transition matrix-based
approach to study  finite-time stabilization of system
(\ref{system3}).  Below, we focus our attention on Lyapunov
function-based approach to present some sufficient   criteria for
finite-time stability and stabilization that are easily tested via
Matlab LMI Toolbox.

{\bf Theorem 4}\label{th4} System (\ref{system1})
is finite-time stable with respect to $(c_1, c_2, T, \{R_k\}_{k\in
N_T})$, if  there exist   a scalar $\alpha\ge 0$,
 and  a
symmetric positive definite matrix sequence $\{P_k>0\}_{k\in N_T}$
solving the following matrix inequalities:
\begin{eqnarray}
&&(\alpha+1)^Tc_1\overline{\lambda}-c_2\underline{\lambda}<0,\label{119}\\
&&A_{j}'R_{j+1}^{\frac{1}{2}}P_{j+1}R_{j+1}^{\frac{1}{2}}A_{j}+C_{j}'R_{j+1}^{\frac{1}{2}}P_{j+1}R_{j+1}^{\frac{1}{2}}C_{j}
-(\alpha+1)R_{j}^{\frac{1}{2}}P_{j}R_{j}^{\frac{1}{2}}<0, \ j\in
N_{T-1}, \label{116}
\end{eqnarray}
where $\underline{\lambda}=\min_{k\in {N_T}}\lambda_{\min} (P_k)$,
$\overline{\lambda}=\lambda_{\max} (P_0)$.

{\textbf {Proof.}} Choose a Lyapunov function  for system
(\ref{system1}) as
$$
V(x,k)=x'\tilde{P}_kx,
$$
where $\tilde{P}_k=R_k^{\frac{1}{2}}P_kR_k^{\frac{1}{2}}$.   Let
$\triangle V(x_j,j)=V(x_{j+1},j+1)-V(x_j,j)$, $j\in N_{T-1}$. By
(\ref{116}),
$$
A_j'\tilde{P}_{j+1}A_j+C_j'\tilde{P}_{j+1}C_j-\tilde{P}_j<\alpha
\tilde{P}_j,
$$
which implies
\begin{eqnarray}
&&E\triangle V(x_j,j)=E[V(x_{j+1},j+1)-V(x_j,j)]\nonumber\\
&=&E[x_j'(A_j'\tilde{P}_{j+1}A_j+C_j'\tilde{P}_jC_j-\tilde{P}_j)x_j]
\le  \alpha  EV(x_j,j). \label{eq fdfn}
\end{eqnarray}
By  (\ref{eq fdfn}), it yields that
\begin{eqnarray}
E[V(x_{j+1},j+1)]&\le& (\alpha+1)  EV(x_j,j) < \cdots\nonumber\\
& \le & (\alpha+1)^{j+1}V(x_0,0)\nonumber\\
&\leq & (\alpha+1)^{T}x_0'\tilde{P}_0x_0. \label{eq werrf}
\end{eqnarray}
It
follows  that
\begin{eqnarray}
EV(x_k,k)&=&E[x_k'\tilde{P}_kx_k]=E\{x'_kR_k^{\frac{1}{2}}P_k R_k^{\frac{1}{2}}x_k\}
\geq \underline{\lambda} E[x'_kR_kx_k] \label{eq cgcvg}
\end{eqnarray}
and
\begin{eqnarray}
V(x_0,0)&=&x_0'\tilde{P}_0x_0=\{x'_0R_0^{\frac{1}{2}}P_0 R_0^{\frac{1}{2}}x_0\}
\leq \overline{\lambda}x'_0R_0x_0. \label{eq fdef}
\end{eqnarray}
So from (\ref{eq werrf})-(\ref{eq fdef}), we have
\begin{eqnarray}
\underline{\lambda}E[x'_kR_kx_k]&\le&(\alpha+1)^k\overline{\lambda}x'_0R_0x_0\nonumber\\
&\leq&(\alpha+1)^k\overline{\lambda}c_1.\nonumber\\
&\leq& (\alpha+1)^T\overline{\lambda}c_1. \label{eq fdvfhj}
\end{eqnarray}
By  (\ref{119}), $E[x'_kR_kx_k]<c_2$ for $k\in N_T$,  which means
that system (\ref{system1}) is finite-time stable with respect to
$(c_1, c_2, T, \{R_k\}_{k\in N_T})$. This theorem is shown.
$\square$

The inequalities (\ref{119})-(\ref{116}) are not LMIs due to
the appearances of $(\alpha+1)^T\overline{\lambda}$ and
$(\alpha+1)P_{j}$, $j\in N_{T-1}$, which leads to a hard computation
in using LMI Toolbox to  solve   (\ref{119})-(\ref{116}).
However, when $\alpha=0$, the inequalities (\ref{119})-(\ref{116}) becomes LMIs.

{\bf Corollary 2}\label{cor4.1}
System (\ref{system1}) is finite-time stable with respect to $(c_1,
c_2, T, \{R_k\}_{k\in N_T})$, if  there exists  a symmetric positive
definite matrix sequence $\{P_k>0\}_{k\in N_T}$ solving the
following LMIs:
\begin{eqnarray}
&& c_1\overline{\lambda}-c_2\underline{\lambda}<0,\label{119cor}\\
&&A_{j}'R_{j+1}^{\frac{1}{2}}P_{j+1}R_{j+1}^{\frac{1}{2}}A_{j}+C_{j}'R_{j+1}^{\frac{1}{2}}P_{j+1}R_{j+1}^{\frac{1}{2}}C_{j} -R_{j}^{\frac{1}{2}}P_{j}R_{j}^{\frac{1}{2}}<0,
\ j\in N_{T-1}, \label{116de}
\end{eqnarray}
where $\underline{\lambda}=\min_{k\in {N_T}}\lambda_{\min} (P_k)$,
$\overline{\lambda}=\lambda_{\max} (P_0)$.

Repeating the above proof and noting that when $-1<\alpha<0$,
(\ref{eq fdvfhj}) should be modified as
\begin{eqnarray}
\underline{\lambda}E[x'_kR_kx_k]&\leq&(\alpha+1)^k\overline{\lambda}x'_0R_0x_0\nonumber\\
&\leq&(\alpha+1)^k\overline{\lambda}c_1.\nonumber\\
&\leq& (\alpha+1)\overline{\lambda}c_1, \ \forall k\in
\{1,2,\cdots,T\}. \label{eq fdvfhjgr}
\end{eqnarray}
We immediately obtain the following result:
{\bf Theorem}\label{th4nn} System (\ref{system1})
is finite-time stable with respect to $(c_1, c_2, T, \{R_k\}_{k\in
N_T})$, if there exist   a scalar $\alpha$ with $-1<\alpha<0$, and a
symmetric positive definite matrix sequence $\{P_k>0\}_{k\in N_T}$
solving the following matrix inequalities:
\begin{eqnarray}
&&(\alpha+1)c_1\overline{\lambda}-c_2\underline{\lambda}<0,\label{119fdfd}\\
&&A_{j}'R_{j+1}^{\frac{1}{2}}P_{j+1}R_{j+1}^{\frac{1}{2}}A_{j}+C_{j}'R_{j+1}^{\frac{1}{2}}P_{j+1}R_{j+1}^{\frac{1}{2}}C_{j}
-(\alpha+1)R_{j}^{\frac{1}{2}}P_{j}R_{j}^{\frac{1}{2}}<0, \ j\in
N_{T-1}, \label{116fd}
\end{eqnarray}
where $\underline{\lambda}=\min_{k\in {N_T}}\lambda_{\min} (P_k)$,
$\overline{\lambda}=\lambda_{\max} (P_0)$.

{\bf Remark 2}
It is easy to see that when $\alpha\le -1$,  (\ref{116fd}) does not
admits solutions $\{P_k>0\}_{k\in N_T}$.

{\bf Theorem 5}\label{th5}
System (\ref{system3}) is  finite-time stabilizable  with respect to
$(c_1, c_2, T, \{R_k\}_{N_T})$ via a linear  state feedback
$u_k=K_kx_k$, if for all $i\in N_T$ and $j\in N_{T-1}$,  there
exists   a   symmetric matrix  sequence ${X_j}$, ${Y_j}$, scalars
$\alpha\ge 0$, $\hat{\lambda}_1>0$, $\hat{\lambda}_2>0$, and $\hat{\lambda}_1 \geq\hat{\lambda}_2$, such that
\begin{eqnarray}
&&\hat{\lambda}_2R_i^{-1}\leq X_i\leq\hat{\lambda}_1R_i^{-1},\label{120fef}\\
&&(\alpha+1)^Tc_1\hat{\lambda}_1-c_2\hat{\lambda}_2<0,\label{121ferf}\\
&&\left[\begin{array}{cccccc}
-(\alpha+1)X_j&(A_jX_j+B_jY_j)'&(C_jX_j+D_jY_j)'\\
*&-X_{j+1}&0\\
*&*&-X_{j+1}
\end{array}\right]
<0. \label{122aerr}
\end{eqnarray}
In addition, if there is a feasible solution for conditions
(\ref{120fef}-\ref{122aerr}), the controller gain can be computed by
$$
K_j=Y_jX_j^{-1}, j\in N_T.
$$

{\textbf {Proof.}} By Definition~\ref{def2} and Theorem 4,
system (\ref{system3}) is finite-time stabilizable  with respect to
$(c_1, c_2, T, \{R_k\}_{N_T})$ via  linear  state feedback
controllers $\{u_k=K_kx_k\}_{k\in N_{T-1}}$ if the following matrix
inequalities are admissible with respect to $\{\alpha\ge 0$,
$\overline{\lambda}>0, \underline{\lambda}>0, P_k>0, k\in {N_T}\}$.
\begin{eqnarray}
&& \underline{\lambda}I<P_i<\overline{\lambda}I,\ i\in {N_T},\label{eq cdcfd}\\
&&(\alpha+1)^Tc_1\overline{\lambda}-c_2\underline{\lambda}<0,\label{119lkp}\\
&&(A_{j}+B_jK_j)'\tilde{P}_{j+1}(A_{j}+B_jK_j) +(C_{j}+D_jK_j)'\tilde{P}_{j+1}(C_{j}+D_jK_j)\nonumber\\
&&\ \ -(\alpha+1)\tilde{P}_j<0, \ j\in N_{T-1}. \label{116v cfd}
\end{eqnarray}
Setting $X_j:=\tilde{P}_j^{-1}$,  $Y_j:=K_jX_j$,  pre- and
post-multiplying  $X_j=X_j'$, the inequality (\ref{116v cfd})
becomes
\begin{eqnarray}\label{118}
&&(A_jX_j+B_jY_j)'X_j^{-1}(A_jX_j+B_jY_j)\nonumber\\
&&+(C_jX_j+D_jY_j)X_j^{-1}(C_jX_j+D_jY_j)-(\alpha+1)X_j<0.
\end{eqnarray}
which is equivalent to (\ref{122aerr}) according to  Schur's
complement. Let $\hat{\lambda}_1=1/{\underline{\lambda}}$,
$\hat{\lambda}_2=1/{\overline{\lambda}}$, and consider
$\tilde{P}_j=R^{1/2}P_jR^{1/2}$,  then the inequalities (\ref{eq
cdcfd}) and (\ref{119lkp}) yield (\ref{120fef}) and (\ref{121ferf}),
respectively. This theorem is proved. $\square$

Unfortunately, the coupled inequalities
(\ref{120fef})-(\ref{122aerr}) are not LMIs.  However, when
$\alpha=0$, Theorem 5 leads to the following easily testing
 finite-time stabilization conditions.
{\bf Corollary}\label{cor4.3}
System (\ref{system3}) is  finite-time stabilizable  with respect to
$(c_1, c_2, T, \{R_k\}_{N_T})$ via a linear  state feedback
$u_k=K_kx_k$, if for all $i\in N_T$ and $j\in N_{T-1}$,  there
exists   a   symmetric matrix  sequence ${X_j}$, ${Y_j}$, positive
scalars
 $\hat{\lambda}_1>0$,  $\hat{\lambda}_2>0$,   $\hat{\lambda}_1 \geq\hat{\lambda}_2$   solving the following
LMIs:
\begin{eqnarray}
&&\hat{\lambda}_2R_i^{-1}\leq X_i\leq\hat{\lambda}_1R_i^{-1},\label{120fhgf}\\
&&c_1\hat{\lambda}_1-c_2\hat{\lambda}_2<0,\label{121fehhhh}\\
&&\left[\begin{array}{cccccc}
-X_j&(A_jX_j+B_jY_j)'&(C_jX_j+D_jY_j)'\\
*&-X_{j+1}&0\\
*&*&-X_{j+1}
\end{array}\right]<0. \label{122aehhhhhr}
\end{eqnarray}
In this case, the desired controller gains is given  by
$$
K_j=Y_jX_j^{-1}, \ j\in N_{T}.
$$

\section{Numerical  Examples\label{sec:EX}}

\hspace{0.13in}  In this section, we give two examples  to show the
effectiveness  of our  proposed results.

{\bf Example 1}
Given the scalars $c_1=0.25$, $c_2=8$, $T=2$, $R_0=1$, $R_1=2$,
$R_2=2$. Consider system (\ref{system1}) with  parameters as
$A_0=C_0=1$, $A_1=C_1=2$ and $x_0=0.5$. By Lemma 1, the
state transition matrices are  computed as
$$
\phi_{1,0}=\left[\begin{array}{ccc}A_0\\C_0\end{array}\right]=
\left[\begin{array}{ccc}1\\1\end{array}\right],
$$
$$
\phi_{2,0}=\left[\begin{array}{ccccccc}
\left[\begin{array}{ccccccc}A_1&0\\0&A_1\end{array}\right]\phi_{1,0}\\
\left[\begin{array}{ccccccc}C_1&0\\0&C_1\end{array}\right]\phi_{1,0}\end{array}\right]
=\left[\begin{array}{ccccccc}2\\2\\2\\2\end{array}\right].
$$
By Theorem 1,
$$
\phi_{0,0}'({R_0})\phi_{0,0}=R_0\leq \frac{c_2}{c_1}R_0,
$$
$$\phi_{1,0}'({I_{2}\otimes R_1})\phi_{1,0}=4\leq \frac{c_2}{c_1}R_0,
$$
$$
\phi_{2,0}'({I_{4}\otimes R_2})\phi_{2,0}=32=\frac{c_2}{c_1}R_0.
$$
Therefore, the given system is  finite-time stable with respect to
$(0.25, 8, 2, \{1,2,2\})$.

\begin{figure}[h]
  \centering
  \includegraphics[width=2.9in]{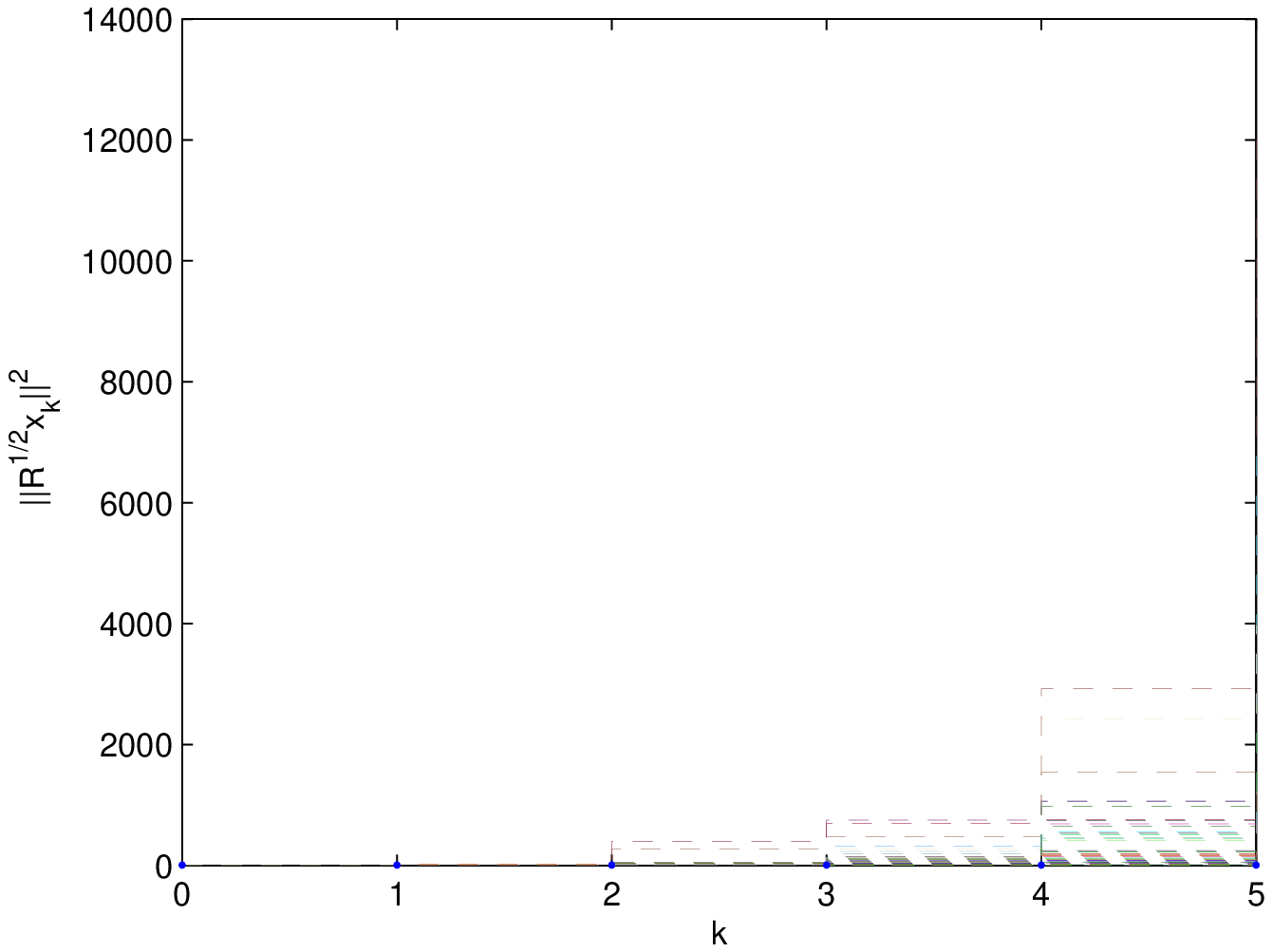}\\[0mm]
  \centering{The response of $R_kx^2_k$ in system (\ref{system1}).}\\[3mm]
  \label{openxrx}
\end{figure}

\begin{figure}[h]
  \centering
  \includegraphics[width=2.9in]{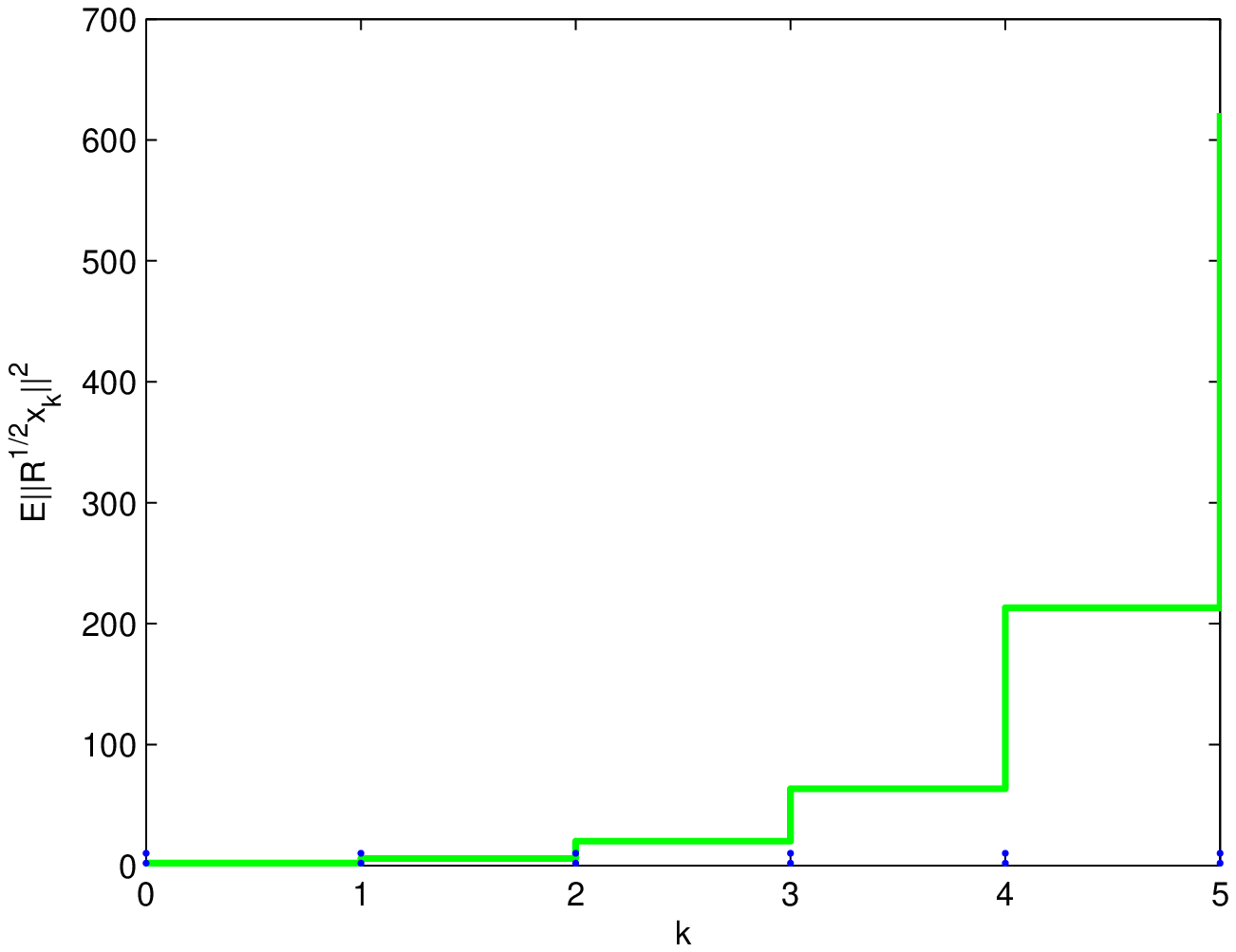}\\[0mm]
  \centering{The response of $E[R_kx^2_k]$ in system (\ref{system1}).}\\[3mm]
  \label{openmean}
\end{figure}

{\bf Example 2}
Given the parameters  $c_1=2, c_2=10, T=20$, $R_k=I_2$, $k\in N_{20}$ with the initial state $x_0=[1 \ 1]'$ in system (\ref{system3}). Assume the  system
(\ref{system3}) is a periodic system with coefficient matrices as
\begin{eqnarray*}
&A_{2i+1}=\left[
\begin{array}{cccc}
-1.023&0.195\\
1.152&0.610
\end{array}
\right],\ \ A_{2i}=\left[
\begin{array}{cccc}
-0.204&1.255\\
0.666&0.282
\end{array}
\right],\\
&C_{2i+1}=\left[
\begin{array}{cccc}
-0.409&1.742\\
-0.482&-0.914
\end{array}
\right],
\ \ C_{2i}=\left[
\begin{array}{cccc}
0.371&0.942\\
0.326&1.748
\end{array}
\right],\\
&B_{2i}=B_{2i+1}=D_{2i}=D_{2i+1}=\left[
\begin{array}{cccc}
1\\
1
\end{array}
\right],
\end{eqnarray*}
where $i\in\{0,1,\cdots,\frac{T}{2}-1=9\}$. We use Matlab to
simulate the system state trajectories 1000 times to obtain  the
approximate value of $E[x_kR_kx_k]$. Figs. \ref{openxrx} and
\ref{openmean} show the  responses of  $x_kR_kx_k$ and
$E[x_kR_kx_k]$  of the uncontrolled  system (\ref{system1}),
respectively. From Figs.  \ref{openxrx} and \ref{openmean}, it can
be seen that the system state is divergent.  It is worth pointing
out that the curves of different colors in Fig. \ref{openxrx}
represent  different experiment results.

\begin{figure}[h]
  \centering
  \includegraphics[width=2.9in]{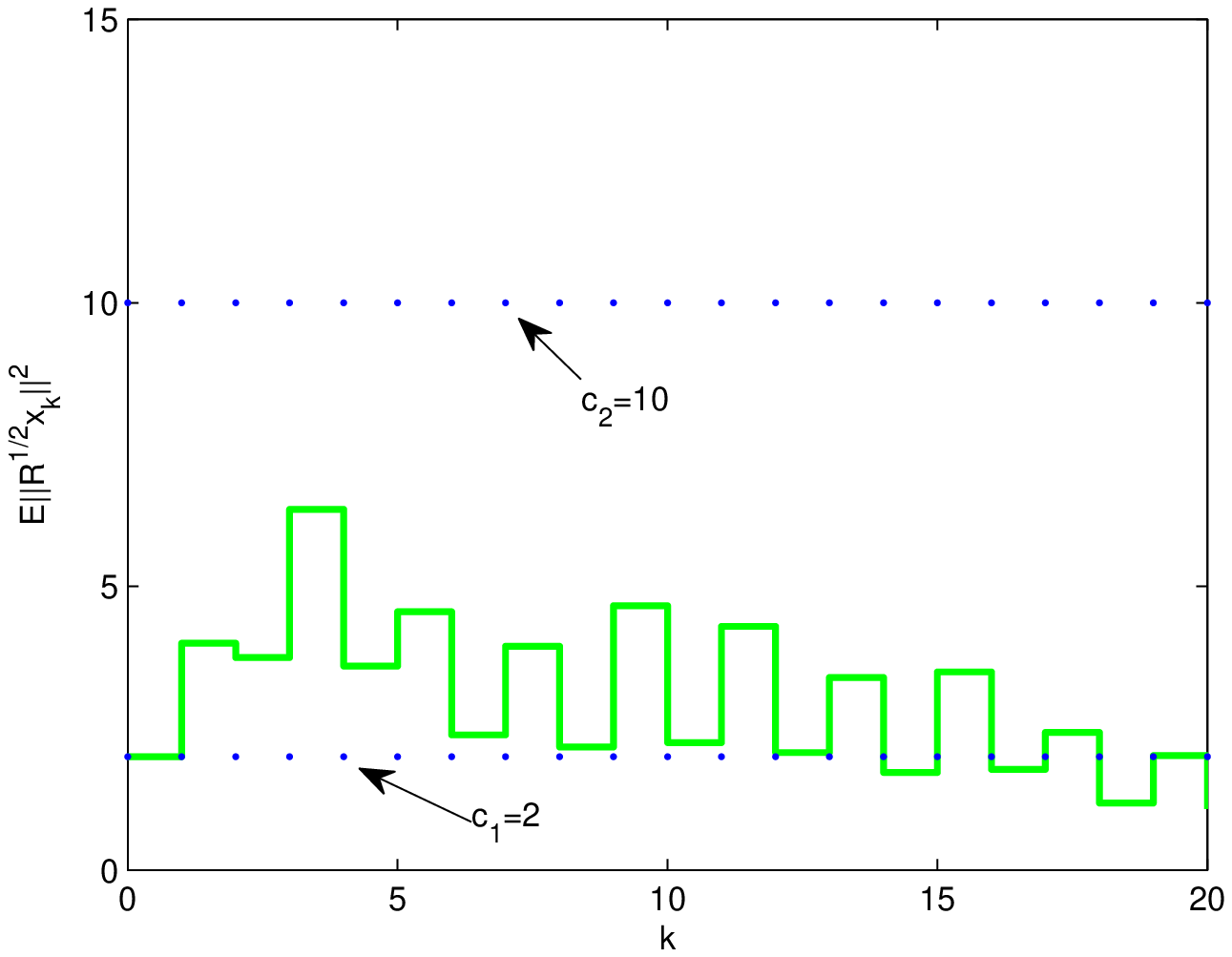}\\[0mm]
  \centering{The response of $E[x'_kR_kx_k]$ in system (\ref{system3}) under  $u_k=K_kx_k$.}\\[3mm]
  \label{mean}
\end{figure}

By Corollary \ref{cor4.3},  we can find a set of feasible solutions to (\ref{120fhgf})-(\ref{122aehhhhhr})  as follows:
$$
X_{2i+1}=\left[
\begin{array}{cccc}
42.710&-17.677\\
-17.677&166.762
\end{array}
\right],
\ \
X_{2i}=\left[
\begin{array}{cccc}
230.590&91.906\\
91.906&170.080
\end{array}
\right],
\ \
Y_{2i+1}=\left[
\begin{array}{cccc}
18.392&-141.356
\end{array}
\right],
$$
$$
Y_{2i}=\left[
\begin{array}{cccc}
-141.855&-206.872
\end{array}
\right],
\ \
\hat{\lambda}_2=0.0033, \ \ \hat{\lambda}_1=0.0249.
$$
Hence, the  periodic feedback gain matrices  are  given by
\begin{eqnarray}
K_k=\left\{\begin{array}{cccc} \left[\begin{array}{ccc} 0.0835&-0.8388\end{array}\right],\ \ \ \ \ \ \ \ \ \ k=2i+1,\\ \\
\left[\begin{array}{ccc}-0.1855&-1.0781\end{array}\right],\ \ \ \ \ \
\ \ \ \ k=2i.
   \end{array}
   \right.
\end{eqnarray}
Under the following state feedback controllers
$$
u_k=K_kx_k, \ k\in N_T,
$$
the closed-loop system of (\ref{system3}) is stabilizable with respect
to $(2, 10, 20, \{R_k=I_2\}_{k\in {N_{20}}})$.
 From Fig.
\ref{mean},   we can see that $E[x'_kR_kx_k]<c_2=10$  for any $k\in
N_{20}$.
\begin{figure}[h]
  \centering
  \includegraphics[width=2.9in]{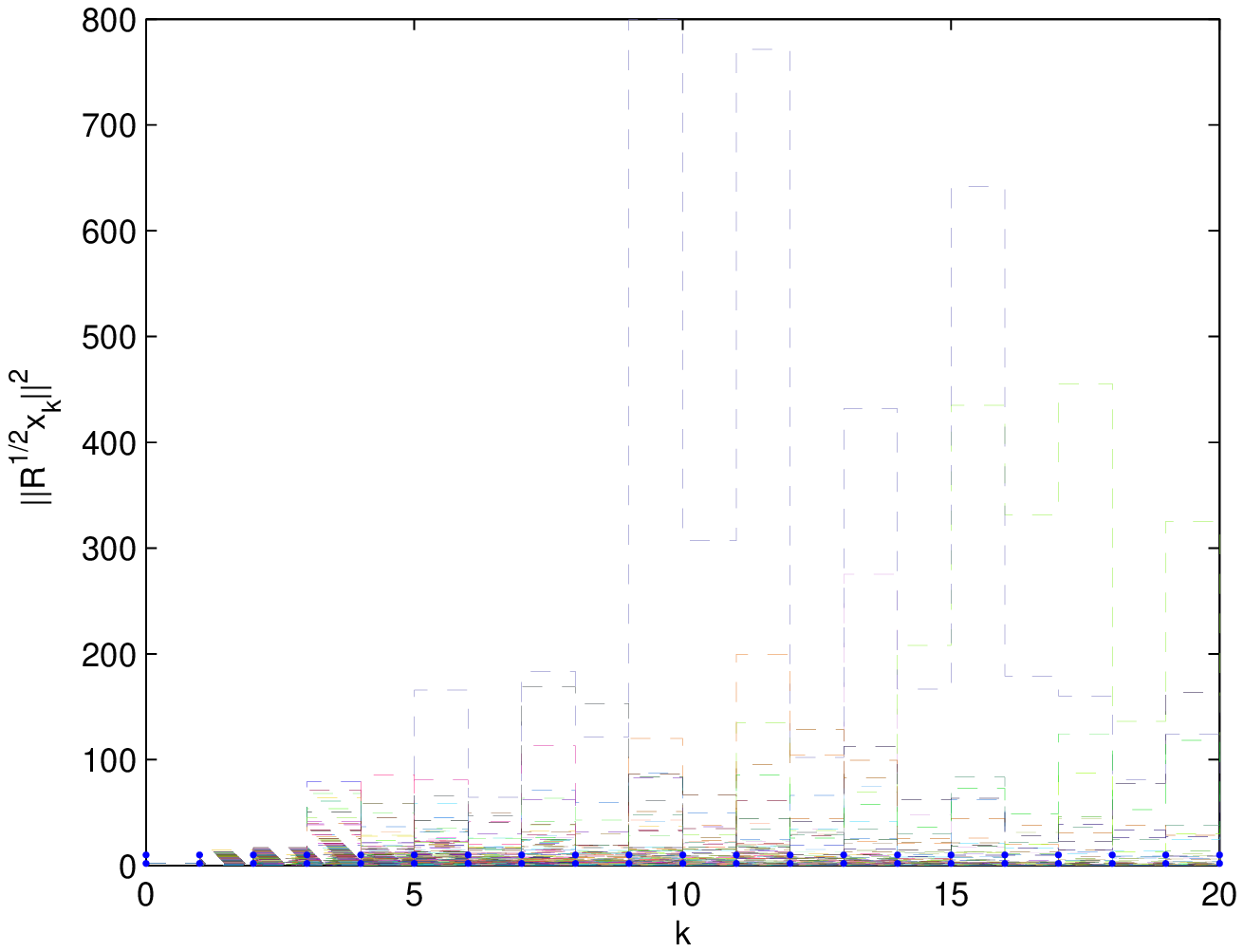}\\[0mm]
  \centering{The response of $x'_kR_kx_k$ in system (\ref{system3}) under  $u_k=K_kx_k$.}\\[3mm]
  \label{xRx}
\end{figure}

%
%

\section{Conclusions\label{sec:Concl}}
\hspace{0.13in}This paper has studied the finite-time stability and
stabilization of linear discrete time-varying stochastic systems.
Necessary and sufficient conditions of finite-time stability and
stabilization have been given based on the state transition matrix.
Meanwhile, sufficient conditions that can be verified using Matlab
LMI Toolbox have also been given. Two examples have  been supplied
to show the effectiveness of our main results.

\end{document}